\documentclass[12pt]{amsart}
\usepackage[centertags]{amsmath}
\usepackage{amsfonts}
\usepackage{amssymb}

\title[Fixed sets]{Fixed sets of automorphisms of   \\ countable, arithmetically saturated structures}

\date{\today}

\author{James H. Schmerl}

\DeclareMathOperator{\aut}{Aut}
\DeclareMathOperator{\fix}{fix}

\DeclareMathOperator{\tp}{tp}

\DeclareMathOperator{\acl}{acl}

\DeclareMathOperator{\dom}{dom}
\DeclareMathOperator{\ran}{ran}

\def\tri{\vartriangleleft}

\newcommand{\ZZ}{{\mathbb Z}}

\def\into{\longrightarrow}

\def\pa{{\sf PA}}
\def\MM{{\mathcal M}}
\def\NN{{\mathcal N}}

\begin{document}

\begin{abstract} If an automorphism $f$  of a structure $\MM$ is such that $\fix(f^k) = \fix(f)$ for all positive $k$, then $\MM | \fix(f)$ is a substructure of $\MM$. The possible isomorphism types of 
$\MM |\fix(f)$ are characterized when $\MM$ is countable and arithmetically saturated. 
\end{abstract}

\maketitle

\bigskip

The purpose of this note is to prove the following theorem, which we are calling the {\it Fix Theorem}.

\bigskip

{\sc   Fix Theorem}: {\em If $\MM$ is a  countable, arithmetically saturated structure  and $D \subseteq M$ is algebraically closed, then there  is $f \in \aut(\MM)$  such that 
$\MM | \fix(f)  \cong \MM | D$ $($even elementarily isomorphic$)$ and  $\fix(f) = \fix(f^k)$ for all $k > 0$.}
 
 \bigskip

 The origin of this theorem can be traced back to  1991 when  Kaye, Kossak \& Kotlarski  \cite[Th.~5.3]{kkk} proved 
 that if $\MM$ is a countable, arithmetically saturated model of  Peano Arithmetic (\pa), then there is  an automorphism $f \in \aut(\MM)$ that moves every undefinable element of~$\MM$ -- that is, $\fix(f) = \acl(\varnothing)$, where $\fix(f) = \{ x \in M : f(x) = x\}$  and 
 $\acl(X)$ is the {\it algebraic closure} of the set $X \subseteq M$. Specifically, for $\MM$  a model  of $\pa$,  $\acl(\varnothing)$ is  the prime elementary substructure  of~$\MM$.  Furthermore, they showed that for any finite $X \subseteq M$, there is $f \in \aut(\MM)$ such that $\fix(f) = \acl(X)$. 
 As a counterpoint to this theorem, Kossak  \cite[Th.~2.6]{ko97} proved that  if $\MM$ is a countable, recursively  saturated model of \pa\ that is not arithmetically saturated, then $\MM | \fix(f) \cong \MM$ for {\it every} $f \in \aut(\MM)$.  It also should be noted  (see \cite[Th.~5.7]{kkk}) that if for such an $\MM$, there is $\NN \prec \MM$ (even $\NN \prec_{\sf end} \MM$) such that $\NN \cong \MM$ and for no $f \in \aut(\MM)$ is $\fix(f) = N$.
 
 Although the proof of the  {K}${\mathbf{^3}}$~theorem  apparently made good use  of \pa,   K\"orner \cite[Lemma~1.3]{korn} eliminated  the need for \pa\    by  proving that 
  {\it every}  countable, arithmetically saturated  structure $\MM$ has an automorphism $f$ 
such that $\fix(f) \subseteq \acl(\varnothing)$. K\"{o}rner calls such an $f$, for any $\MM$,  a {\it maximal} automorphism of $\MM$. K\"orner's theorem was subsequently  improved by Duby \cite[Th.~8]{duby} as follows.

\bigskip

{\sc Duby's Theorem}: {\em  If $\MM$ is countable and arithmetically saturated, then there is $f \in \aut(\MM)$ such that   $\fix(f^k) = \acl(\varnothing)$ whenever $1 \leq k < \omega$.}

\bigskip
   This improvement to  K\"orner's theorem is automatic for models of \pa\ and  even for all linearly ordered $\MM$, but not so for some other structures. Duby refers to such an $f$ as an \mbox{$\omega$-{\it maximal}} automorphism of $\MM$. It is readily seen that Duby's Theorem generalizes to: If $\MM$ is countable and arithmetically saturated and $X \subseteq M$ is finite, then  there is $f \in \aut(\MM)$ such that $\fix(f^k) = \acl(X)$ for all $k > 0$. 

The question was raised  (\cite[Question~2.7]{ko97} and  \cite[Question~9, p.\@ 291]{ksbook}) as to what other possibilities there are for the isomorphism types of $\MM | \fix(f)$ when $\MM$ is  a countable,  arithmetically saturated model of \pa\ and $f \in \aut(\MM)$.  According to Enayat \cite[\S4.2]{en07}, I had earlier conjectured that if  $\MM$ is a countable,  arithmetically saturated model of $\pa$ and $\NN \preccurlyeq \MM$, then there is $f \in \aut(\MM)$ such that $\MM |\fix(f) \cong \NN$. Kossak \cite[Th.~2.8]{ko97} lent credence to this conjecture by observing that for every countable model $\NN$  of $\pa$, there are a countable, arithmetically saturated $\MM \succ \NN$ and  $f \in \aut(\MM)$ such that $\fix(f) = N$. Incidentally,  K\"orner's theorem easily implies that for every countable $\NN$ (not necessarily a model of \pa), there are a countable, arithmetically saturated $\MM \succ \NN$ and $f \in \aut(\MM)$ such that $\fix(f) = N$.{\footnote{Assume $\NN$ is infinite. Let $<$ be a linear order of $N$ having ordertype 
$\omega$, and let $(\MM,<) \succ (\NN,<)$ be countable and arithmetically saturated. Apply K\"orner's theorem  to $(\MM,<)$.}} My conjecture was later confirmed by Enayat \cite[Th.~4.2.1]{en07} 
using iterated ultrapowers.

\bigskip

{\sc Enayat's Theorem}: {\em If $\MM \models \pa$ is a countable and arithmetically saturated and $\NN \preccurlyeq \MM$, then there is $f \in \aut(\MM)$ such that $\MM | \fix(f) \cong \NN$.}

\bigskip

 The Fix Theorem  is a common generalization of Duby's and Enayat's Theorems. 
We note that the Fix Theorem is best possible in the sense that if $\MM$ is  any structure and $f \in \aut(\MM)$ is $\omega$-maximal, then $\MM | \fix(f)$ is algebraically closed.

 Following this introduction there are three sections. The first contains preliminary material, including those definitions  that are  needed to understand the Fix  Theorem.  Most readers can probably skip this section. The second section  contains  proofs of Duby's Theorem and of Enayat's Theorem. The Fix Theorem is proved in  the third section. 
 
 \bigskip


{\bf \S1.~Preliminaries.} In order to avoid any possible snags, 
all structures  considered here are for a finite language -- that is, for each structure $\MM$, there is a finite language~${\mathcal L}$ such that $\MM$ is an ${\mathcal L}$-structure. Typically, an ${\mathcal L}$-structure will be denoted by a (possibly adorned) script letter, such as  $\MM$, and then its universe is understood to be denoted by the corresponding (similarly adorned)  latin letter, such as~$M$. If $A \subseteq M$, then ${\mathcal L}(A)$ is the language ${\mathcal L}$ adjoined with a constant symbol for each $a \in A$. 

As usual, $\omega$ is the set of natural numbers (i.e. nonnegative integers) and ${\mathbb Z}$ is the set of all (negative and nonnegative) integers. Hence, we have that $\omega \subseteq {\mathbb Z}$. We typically write $n < \omega$ instead of $n \in \omega$.

   If $a \in M$ and $X \subseteq M$, then  $\tp( a/ X)$, the {\bf type of} $a$ {\bf over} $X$,  is the set of all $1$-ary ${\mathcal L}(X)$-formulas $\varphi(x)$ such that $\MM \models \varphi( a)$. As usual, $\tp(a) = \tp(a / \varnothing)$. 
   
Let $\MM$ be an ${\mathcal L}$-structure. The set (or group) of automorphisms  of  $\MM$ is  $\aut(\MM)$. For $f \in \aut(\MM)$, its {\bf fixed set} is $\fix(f) = \{x \in M : f(x) = x\}$.   If $X \subseteq M$, then the {\bf algebraic closure} of $X$, denoted by $\acl(X)$, is the union of all finite $D$ that are definable in $\MM$ by an ${\mathcal L}(X)$-formula.  Easily, $\acl\big(\acl(X)\big) = \acl(X)$. We write $\acl(a)$ instead of $\acl(\{a\})$. A subset $X \subseteq M$ is {\bf algebraically closed} if  $X = \acl(X)$.    If $A,X \subseteq M$ and $\sigma : A \into M$, then $\sigma$ is {\bf elementary over} $X$ if, whenever $\varphi(\overline x)$ is an $n$-ary ${\mathcal L}(X)$-formula and $\overline a \in A^n$, then $\MM \models \varphi(\overline a) \leftrightarrow \varphi\big(\sigma(\overline a)\big)$.    If $D_1,D_2 \subseteq M$, then $D_1$ is {\bf elementarily isomorphic} to $D_2$ if there is an elementary surjection $\sigma : D_1 \into D_2$. 
 The parenthetical part in  the conclusion of the Fix Theorem says that $\fix(f)$ and $D$ are elementarily isomorphic.

  Recall that a structure $\MM$ is {\bf recursively saturated} if whenever $X \subseteq M$ is finite, $\Sigma(x)$ is a computable  (i.e.~recursive) set of 1-ary ${\mathcal L}(X)$-formulas that is finitely realizable  in $\MM$, then $\Sigma(x)$ is realizable in $\MM$. Alternatively, $\MM$ is recursively saturated iff whenever $X \subseteq M$ is finite, $\Sigma(x)$ is a set  of 1-ary ${\mathcal L}(X)$-formulas that is finitely realizable in $\MM$ and is computable  in  some $\tp(\overline a)$, where $\overline a \in M^n$ and $n < \omega$, then $\Sigma(x)$ is realized in~$\MM$. Analogous with this latter characterization of recursive saturation,  $\MM$ is {\bf arithmetically saturated} if whenever $X \subseteq M$ is finite, $\Sigma(x)$ is a set  of 1-ary ${\mathcal L}(X)$-formulas that is finitely realizable in $\MM$ and is arithmetic in some $\tp(\overline a)$, where $\overline a \in M^n$, then  $\Sigma(x)$ is  realized in $\MM$. 
A structure $\MM$ is {\bf resplendent} if whenever $R$ is a new $k$-ary relation symbol, $A \subseteq M$ is finite and   $\sigma$ is an $({\mathcal L}(A) \cup \{R\})$-sentence  that is modeled by an expansion of some $\NN \succcurlyeq \MM$,  then $\MM$ has an expansion modeling  $\sigma$. Every resplendent structure has the following stronger property: Whenever $\Sigma$ is a set of 
$({\mathcal L}(A) \cup \{R\})$-sentences that is computable in $\tp(\overline a)$ for some $\overline a$ and that is modeled by an expansion of some $\NN \succcurlyeq \MM$, then $\MM$ has an expansion modeling $\Sigma$. Every resplendent  structure is recursively saturated. Conversely,  every countable, recursively saturated structure is resplendent. Moreover, every countable, recursively saturated structure is {\bf chronically} resplendent - meaning that it has such an expansion to a recursively saturated model of $\Sigma$. 
We say that $\MM$ is {\bf arithmetically} resplendent if whenever $\Sigma$ is a set of 
$({\mathcal L}(A) \cup \{R\})$-sentences that is arithmetic in some $\tp(\overline a)$ and that is modeled by  some expansion of some $\NN \succcurlyeq \MM$, then $\MM$ has an expansion modeling $\Sigma$.  We  say that $\MM$ is {\bf chronically} arithmetically resplendent if, furthermore, there is such an arithmetically resplendent expansion of $\MM$.    We will  make use of the following well-known lemma.

\bigskip

{\sc Lemma 1.1}: {\em Every countable, arithmetically saturated structure  is chronically arithmetically resplendent.}

\bigskip

 {\bf \S2.~Proving Duby's  and Enayat's Theorems.} In this section   we give   proofs of Duby's Theorem and Enayat's Theorem. The proof of Duby's Theorem  is essentially the one in \cite[\S2]{duby} but perhaps simpler   in its details.   Our proof of Enayat's Theorem is based on our presented   proof of Duby's Theorem and appears to be different from the one in~\cite{en07}. 
 
 Throughout this section, $\MM$ is always a countable,  arithmetically saturated ${\mathcal L}$-structure and  ${\mathcal L}$ is finite.

Duby uses Lemma~2.2 in his proof. To prove that lemma, Duby invokes the following lemma.

\bigskip

{\sc Lemma 2.1}: {\em If $K,X \subseteq M$ are finite and $a \in M$, then there is $d \in M$ such that $\tp\big(d / \acl(K)\big)
= \tp\big(a / \acl(K)\big)$ and 
\begin{equation} \tag{$*$}
\acl(K \cup X) \cap \acl(K \cup \{d\}) = \acl(K).
\end{equation}}

\bigskip

We indicate how Duby proves this lemma in \cite{duby}. Notice that by replacing $\MM$ with $(M,x)_{x\in K}$, which is also countable and arithmetically saturated, we can assume that $K = \varnothing$.{{\footnote{This is done mostly as a matter of notational simplicity. We call this procedure the {\em nullification of }$K$ and will use it often.} Duby's proof of this lemma  in \cite{duby} is in two parts, each involving a set $\Phi(x)$ of ${\mathcal L}(X)$-formulas.  In the first part, he invokes the Separation Theorem of P.~M.~Neumann (\cite[Lemma 1]{duby}) to show in Lemma~3 that $\Phi(x)$ is finitely consistent. In the second part, he shows in \cite[Theorem~6]{duby} that $\Phi(x)$ is arithmetic in the type of $\tp\big(a~ / \acl(X)\big)$.
However, in $\S3$ we will prove Lemma~3.1 from which  Corollary~3.2  and then Lemma~2.1 easily follow.

\bigskip

{\sc Lemma 2.2}: {\em Suppose that $A \subseteq M$ is finite and $\sigma : A \into M$ is elementary over $\acl(\varnothing)$. If $c \in M$, then there is $d \in M$ such that $\sigma \cup \{\langle c,d \rangle\}$ is elementary over $\acl(\varnothing)$ and 
\begin{equation} \tag{$**$}
\acl(A \cup \{c\}) \cap \acl(\sigma[A] \cup \{d\}) \subseteq  \acl(\sigma[A]).
\end{equation}}

\bigskip

{\it Proof}. Let $\sigma, A$ and $c$ be as given. Using the arithmetic saturation of $\MM$, let $a \in M$ be such that $\sigma \cup \{\langle c,a \rangle\}$ is elementary over $\acl(\varnothing)$. Letting  $X = A \cup \{c\}$ and $K = \sigma[A]$, apply Lemma~2.1 to get $d \in M$.
One easily verifies that $d$ is as in ($**$). \qed

\bigskip

We next make a definition of a concept that is perhaps implicit in \cite[\S2]{duby}. To acknowledge the importance of Duby's contribution, we name it after him.

\bigskip

{\sc Definition 2.3}:  We say that $\sigma$ is {\bf Duby} if there is a finite $A \subseteq M$ such that $\sigma : A \into M$ and $(1)$ and (2) hold, where:

(1) $\sigma$ is elementary over $\acl(\varnothing)$;

(2) whenever $\sigma \subseteq f \in \aut(\MM)$ and $\acl(\varnothing) \subseteq \fix(f)$,  then
$$
\bigcap_{k \in \ZZ}\acl(f^k[A]) = \acl(\varnothing).
$$

\bigskip

One easily sees that the empty function $\varnothing$ is Duby. Also, every Duby function is 
one-to-one (by (1)) and  its inverse  is also Duby.

\bigskip

{\sc Lemma 2.4}: {\em Suppose that  $\sigma$ is Duby and that $c \in M$.}

(a)  {\em There is $d \in M$ such that $\sigma \cup \{\langle c,d \rangle\}$ is Duby.} 

(b) {\em There is $d \in M$ such that $\sigma \cup \{\langle d,c \rangle\}$ is Duby.}

\bigskip

{\it Proof}. Let $\sigma : A \into M$ be Duby and $c \in M$. As we noted after Definition 2.3, each Duby function is one-to-one and its inverse  is also Duby. Hence, conclusion (b) follows from (a), so  we prove only (a).

 Let $d$ be as in Lemma~2.2. We wish to show that  $\sigma \cup \{\langle c,d \rangle\}$ is Duby; that is, we want to prove $(1')$ and $(2')$, which are (1) and (2) of Definition~2.3 but with $\sigma \cup \{\langle c,d \rangle\}$ and $A \cup \{c\}$ replacing $\sigma$ and $A$, respectively. It is given that $(1')$ holds. In order to prove $(2')$,  suppose that $\sigma \cup \{\langle c,d \rangle\}  \subseteq f \in \aut(\MM)$. Obviously, 
 $$
 \bigcap_{k\in \ZZ}\acl(f^k[A \cup \{c\}] \supseteq \acl(\varnothing).
 $$
 In proving the reverse inclusion, we will use that for any $X \subseteq M$ and $g \in \aut(\MM)$, then 
 $\acl(g[X]) = g[\acl(X)] \mbox{\ and \ } g[X \cap Y] = g[X] \cap g[Y]$.
 Thus,

\begin{align*}
\bigcap_{k \in \ZZ}\acl(f^k[A \cup \{c\}])&= \bigcap_{k \in \ZZ}\acl(f^k[A \cup \{c\}]) \cap \bigcap_{k \in Z}\acl(f^{k+1}[A \cup  \{c\}]) \\ 
&= \bigcap_{k \in \ZZ}\acl(f^k[A \cup \{c\}]) \cap \acl(f^{k+1}[A \cup  \{c\}]) \\  
& =  \bigcap_{k \in \ZZ} f^k[\acl(A \cup \{c\}) \cap \acl(f[A] \cup \{d\})] \\ 
& =  \bigcap_{k \in \ZZ} f^k[\big(\hspace{-2pt}\acl(A \cup \{c\})\big) \cap \big(\hspace{-2pt}\acl(\sigma[A] \cup \{d\})\big)] \\ 
& \subseteq \bigcap_{k \in \ZZ} f^k[\acl(\sigma[A])] 
 = \bigcap_{k \in \ZZ}\acl(f^k[A])     = \acl(\varnothing). 
\end{align*}
The inclusion in the fifth line is due to  $(**)$ of Lemma~2.2.
Thus, we have that $\sigma \cup \{\langle c,d \rangle\}$ is Duby, completing the proof (a) and, therefore, also of (b). \qed

\bigskip

{\it Proof of Duby's Theorem}. (Refer to the Introduction for a statement of Duby's Theorem.)  Using the countability of $\MM$ and by repeated applications of Lemma~2.4,   we  get an increasing sequence $\sigma_0 \subseteq \sigma_1 \subseteq  \sigma_2 \subseteq \cdots$  of Duby functions such that $\sigma_0 = \varnothing$ and 
 for every $c \in M$, there are $i,j < \omega$ such that $c \in \dom(\sigma_i)$ 
 (by Lemma~2.4(a)) and $c \in \ran(\sigma_j)$ (by Lemma~2.4(b)). We claim that $f = \bigcup_{i<\omega}\sigma_i$ is an $\omega$-maximal automorphism of $\MM$. It follows from Definition~2.3(1) that  $f \in \aut(\MM)$ and that $\fix(f) \supseteq \acl(\varnothing)$.
 
 We prove that $f$ is $\omega$-maximal. Suppose that there are $a \in M$ and positive $k < \omega$ such that $f^k(a) = a$. Let $A = \{f^i(a) : i < k\}$. Let $j$ be such that $\sigma_j : A_j \into M$ and $A \subseteq \dom(\sigma_j) \cap \ran(\sigma_j)$. Then, $\sigma_j \subseteq f \in \aut(\MM)$ and 
 $$
 a \in A \subseteq  \bigcap_{k \in \ZZ}\acl(f^k[A_j]) = \acl(\varnothing).
 $$
Thus, $a \in \acl(\varnothing)$, implying that $f(a)  = a$. This completes the proof of Duby's Theorem.
\qed

\bigskip

We next consider relativizations of Duby's Theorem.

\bigskip

{\sc Corollary~2.5}: {\em If $K \subseteq M$ is finite, then there is $f \in \aut(\MM)$ such that $\fix(f^k) = \acl(K)$ whenever $1 \leq k < \omega$.}

\bigskip
{\it Proof}. The simplest way to prove this is by the nullification of $K$, as  in the proof of  Lemma~2.1. Thus, assume that $K = \varnothing$, and then apply Duby's Theorem. \qed

\bigskip

There is another approach to proving Corollary~2.5.  To see this, we generalize Definition~2.3 as follows.

\bigskip

{\sc Definition 2.6}:  We say that $\sigma$ is {\bf Duby over} $K$ if there is a finite $A \subseteq M$ such that $\sigma : A \into M$ and $(0)$ -- $(2)$ hold, where:

\begin{itemize}

\item[(0)] $K \subseteq M$ and $K$ is finite;

\item[(1)] $\sigma$ is elementary over $\acl(K)$;

\item[(2)] whenever $\sigma \subseteq f \in \aut(\MM)$ and $\acl(K) \subseteq \fix(f)$, then
$$
\bigcap_{k \in \ZZ}\acl(f^k[A \cup K]) = \acl(K).
$$
\end{itemize}

\bigskip

If $K = \varnothing$, then (1) and (2) of the previous definition are the same as (1) and (2) of Definition~2.3. Thus, $\sigma$ is Duby over $\varnothing$ iff $\sigma$ is Duby. It is obvious that  for each finite $K \subseteq M$, $\varnothing$ is Duby over $K$.

The following ``over K" version of Lemma~2.4 can also be used to prove Corollary~2.5.
\bigskip

\bigskip

{\sc Lemma 2.7}: {\em Suppose that  $\sigma$ is Duby over $K$  and that $c \in M$.}

(a)  {\em There is $d \in M$ such that $\sigma \cup \{\langle c,d \rangle\}$ is Duby over $K$.} 

(b) {\em There is $d \in M$ such that $\sigma \cup \{\langle d,c \rangle\}$ is Duby over $K$.}
\qed

\bigskip

We now turn to the proof of Enayat's Theorem, for which the following lemma is used.

\bigskip

{\sc Lemma 2.8}: {\em Suppose that $\MM \models \pa$ and that $\sigma$ is Duby over $K$. If $a \in M$, then there is $b \in M$ such that $\tp\big(b / \acl(K)\big) = \tp\big(a / \acl(K)\big)$ and  $\sigma$ is Duby over $K \cup \{b\}$.}

\bigskip

{\it Proof.} Let $\MM,\sigma, K$ and $a$ be as given. 
Let $\varphi(x) \in \tp\big(a~/ \acl(K)\big)$. Then $\MM \models \varphi(a)$, so that we can let $c$ be the least such that $\MM \models \varphi(c)$. Thus, $c \in \acl(K)$, so that
$\sigma$ is Duby over $K \cup \{c\}$. By the arithmetic saturation of $\MM$, there is $b \in M$ such that $\tp\big(b / \acl(K)\big) = \tp\big(a / \acl(K)\big)$ and $\sigma$ is Duby over $K \cup \{b\}$. \qed

\bigskip

{\it Proof of Enayat's Theorem}. (Refer to the Introduction for a statement of Enayat's Theorem.) For this proof, $\MM$ is a model of \pa, which, as usual, is countable and arithmetically saturated. Obviously, both $\MM$ and $\NN$ are  infinite.  Let $a_0,a_1,a_2, \ldots$ and $c_0,c_1,c_2, \ldots$ be non-repeating enumerations of $N$ and $M$, respectively.
We will obtain two sequences: $\langle \sigma_n : n < \omega \rangle$ and $\langle K_n : n < \omega \rangle$. These sequences will be such that whenever $n < \omega$, then:

\begin{itemize}
\item[(1.n)] $\sigma_n$ is Duby over $K_n$;
\item[(2.n)] $\sigma_n \subseteq \sigma_{n+1}$ and $K_n \subseteq K_{n+1}$.

\end{itemize}

We start by letting $\sigma_0 = K_0 = \varnothing$. Obviously, (1.0) is satisfied. We proceed recursively. So, fix some $n < \omega$ and assume that we already have $\sigma_n$ and  $K_n$ satisfying $(1.n)$. Using Lemmas~2.7(a), 2.7(b) and ~2.8, we get $d_n,d_n'$ and $b_n$ such that:

\begin{itemize} 
\item[(3.n)] $\tp\big(b_n~ / \acl(K_n)\big) = \tp\big(a_n~ / \acl(K_n)\big);$
\item[(4.n)] $\sigma_n \cup \{\langle c_n,d_n \rangle, \langle d'_n,c_n \rangle\}$ is Duby over $K_n \cup \{b_n\}$.
\end{itemize}
We then let $\sigma_{n+1} = \sigma_n \cup \{\langle c_n,d_n \rangle, \langle d'_n,c_n \rangle\}$
and $K_{n+1} = K_n \cup \{b_n\}$, so that both (3.n+1) and (4.n) are satisfied. Letting $f = \bigcup_{n<\omega} \sigma_n$, we see, just as in the proof of Duby's Theorem,  that $f \in \aut(\MM)$ and that $\MM | \fix(f) \cong \NN$, as demonstrated by $\{\langle a_n,b_n \rangle : n < \omega\}$.  \qed

  \bigskip

%
%

{\bf \S3.~Proving the Fix Theorem.}  In this section, as in the previous section, $\MM$ is always a countable,  arithmetically saturated \mbox {${\mathcal L}$-structure} and  ${\mathcal L}$ is finite. 

We begin this section with a possible, yet incorrect,  approach to proving the Fix Theorem. Suppose that one could remove from Lemma~2.8 the requirement that $\MM$ be a model of \pa. Then we would be able to prove the Fix Theorem just as we proved Enayat's Theorem. But that does not work, as the following example shows. Let $\MM = (M,R)$, where $R$ is an equivalence relation on $M$ having exactly two  equivalence classes, both of which are infinite. Notice that $\MM$ is $\aleph_0$-categorical and, hence, arithmetically saturated. Let $c,d \in M$ be in different equivalence classes. Then $\sigma = \{\langle c,d \rangle\}$ is Duby over $\varnothing$, but there is no $b \in M$ such that 
$\tp(c,b) = \tp(d,b)$.

However, what we will do to prove the Fix Theorem is to strengthen the definition of {\it Duby over } $K$ in Definition 2.6 to {\it superDuby over } $K$ in Definition~3.10. We then will show that there are three lemmas concerning this stronger notion, namely Lemmas~3.11 -- 3.13, from which  we can prove the Fix Theorem in a way that is very similar to the given proof of Enayat's Theorem. 
 
Recall that we had promised  Lemma~3.1 from which  Lemma~2.1 easily follows.

\bigskip

{\sc Lemma 3.1}: {\em Suppose that $K,X \subseteq M$ are finite. Then there is $\NN \preccurlyeq \MM$ such that  
\begin{equation} \tag{$***$}
\acl(K \cup X) \cap N =  \acl(K)
\end{equation}
and such that $(\MM,N)$  is arithmetically saturated.}

\bigskip

{\it Proof}. Let $K$ and $X$ be as given.  If $\MM$ is finite, then $\acl(\varnothing) = M$, so we can let $\NN = \MM$.  Hence, we assume that $\MM$ is infinite. 
By the nullification of $K$, we also assume that $K = \varnothing$. 

Let $\Gamma_0$ be a computable set of 
$({\mathcal L} \cup \{N\})$-sentences asserting that $\NN \prec \MM$.  Let $\Gamma_1$ be the set of $({\mathcal L}(X) \cup \{N\})$-sentences consisting of just those sentences
$$
\forall x[\theta_1(x) \rightarrow x \not\in N],
$$
where $\theta_1(x)$ is an algebraic ${\mathcal L}(X)$-formula such that
$$
\MM \models \forall x \big(\theta_1(x) \rightarrow \neg\theta_0(x)\big)
$$
for every algebraic ${\mathcal L}$-formula $\theta_0(x)$. The set $\Gamma_1$ is arithmetic in the set of true 
${\mathcal L}(X)$-sentences. Let $\Gamma = \Gamma_0 \cup \Gamma_1$. If $N \subseteq M$, $\NN = \MM|N$ and $(\MM,N) \models \Gamma$, then $\NN \preccurlyeq \MM$ and 
$\acl(X) \cap N = \acl(\varnothing)$. Thus, we want such an $\NN$ with $(\MM,N)$ being arithmetically saturated. 

We will show that there is $N \subseteq M$ such that $(\MM,N) \models \Gamma$.
Let $(\MM_0,N_0) \succ (\MM,N)$ be uncountable and saturated, and let $\NN_1 \prec \NN_0$ be countable.  Let ${\mathbf T}$ be an uncountable set of elementary substructures of $\MM_0$ such that $\NN_1 \in {\bf T}$ and whenever $\NN_2,\NN_3 \in {\mathbf T}$ are distinct, then $N_2 \cap N_3 = \acl(\varnothing)$ and there is $\alpha \in \aut(\MM_0,N_0)$ such that $N_2 = \alpha[N_3]$. Since ${\mathbf T}$ is uncountable and $\acl(X)$ is countable, there is $\NN_2 \in {\mathbf T}$ such that $N_2 \cap \acl(X) = \acl(\varnothing)$. So there is $\NN_2 \prec \MM_0$ such that $(\MM_0,N_2) \models \Gamma$. 

Then,  there is $\NN \prec \MM$ such that $(\MM,N) \models \Gamma$. But then there is such an $\NN$ which, in addition, by Lemma~1.1, is such that $(\MM,N)$ is arithmetically saturated.
 \qed
 
 \bigskip
 
 Throughout this section, we will use $(*)$, $(**)$ and $(***)$ just as they are in Lemmas~2.1, 2.2 and~3.1, respectively.

\bigskip

{\sc Corollary 3.2}: {\em If $K,X \subseteq M$ are finite, then there is $\NN \preccurlyeq \MM$ such that  for each $d \in N$,
\begin{equation} \tag{$*$}
\acl(K \cup X) \cap \acl(K \cup \{d\}) = \acl(K)
\end{equation}
and such that $(\MM,N)$ is arithmetically saturated.}

\bigskip

{\it Proof}. Let $K,X \subseteq M$  be finite.   Lemma~3.1 asserts the existence of  $\NN \preccurlyeq \MM$ such that $\acl(K \cup X) \cap N = \acl(K)$ and  $(\MM,N)$ is arithmetically saturated.  Then every $d \in N$ satisfies $(*)$.
 \qed
 
 \bigskip
 
 {\it Proof of Lemma~2.1}.  In Corollary~3.2, let $K = \varnothing$ and then let $d \in N$ be such that $\tp\big(d / \acl(\varnothing)\big)
= \tp\big(a / \acl(\varnothing)\big)$.  \qed

\bigskip

Lemma~3.1 implies the following strengthening of itself. 

\bigskip

{\sc Corollary 3.3}: {\em Suppose that $\NN' \preccurlyeq \MM$ is such that $(\MM,N')$ is arithmetically saturated. Further, suppose that $K,X \subseteq M$ are finite and that $K \subseteq N'$. Then there is $\NN \preccurlyeq \NN'$ such that 
\begin{equation} \tag{$***$}
\acl(K \cup X) \cap N = \acl(K)
\end{equation}
and such that $(\MM,N)$ is arithmetically saturated.}

\bigskip

{\it Proof}: Let $\NN'$, $K$ and $X$ be as given.  By the nullification of $K$, we assume that $K = \varnothing$. Since $(\MM,N')$ is arithmetically saturated, we can apply Lemma~3.1 to get $(\MM',N) \preccurlyeq (\MM,N')$ such that $(\MM,N',M',N)$ is arithmetically saturated and $\acl(X) \cap M' = \acl(\varnothing)$. But then, $(\MM,N)$ is arithmetically saturated and $\acl(X) \cap N = \acl(\varnothing)$. \qed

\bigskip

There is also the corresponding strengthening of Corollary~3.2. 

\bigskip

{\sc Corollary 3.4}: {\em Suppose that $\NN' \preccurlyeq \MM$ is such that $(\MM,N')$ is arithmetically saturated. Further, suppose that $K,X \subseteq M$ are finite and that $K \subseteq N'$. Then there is $\NN \preccurlyeq \NN'$ such that  for each $d \in N$,
\begin{equation} \tag{$*$}
\acl(K \cup X) \cap \acl(K \cup \{d\}) = \acl(K)
\end{equation}
and such that $(\MM,N)$ is arithmetically saturated.} \qed

\bigskip

Recall that $(T,<)$ is a partially ordered set (or simply a {\it poset}) if $(T,<)$ is such that $<$ is a transitive and irreflexive binary relation on $T$. (In a poset $(T,<)$, we allow the possibility that $T = \varnothing$.) If $(T,<)$ is a poset and $s,t \in T$ are such that $s < t$, but for no $r \in T$ is $s < r < t$, then $t$ is an 
{\it immediate successor} of $s$. Also, $\tri$ is a linear order of $T$ if  $(T,\tri)$ is a poset and whenever $s,t \in T$ are distinct, then either $s \tri t$ or $t \tri s$. If $(T,<)$, $(T, \tri)$ are posets, then $(T,\tri)$ extends $(T,<)$ if $s <t$ implies $s \tri t$. If $(T,<)$ is a poset and  $t \in T$, then we let $D(t) = \{s \in T : s \leq t\}$. 
If $(T,<)$ is a poset, then 
$$
\min(T) = \{t \in T : D(t)  = \{t\}\}.
$$
 We say that the poset $(T,<)$ is a  {\it forest} if  for each $t \in T$, the set $D(t)$ is linearly ordered by $<$. If $(T,<)$ is a finite forest, then a {\bf twig} is a subset $W \subseteq T$ such that either $W = \varnothing$ or $W = D(t)$ for some $t \in T$. 

\bigskip

{\sc Lemma 3.5}: {\em Suppose that $(T,<)$ is a forest and $|T| = n < \omega$. Then there is a linear order $(T, \tri)$ extending $(T,<)$ such that $T = \{t_i : i < n\}$, $t_0 \tri t_1 \tri \cdots \tri t_{n-1}$ and $t_i \tri t_j$ iff $i < j < n$. Furthermore,  
 whenever $i+1 \leq j < n$, then
$$
D(t_i) \cap D(t_j) \subseteq D(t_i) \cap D(t_{i+1}).
$$}

\bigskip

{\it Proof}. Let $|T| = n < \omega$. If $n \leq 1$, then let $\tri$  be $\varnothing$. So, we can assume that $n \geq 2$. Let $t_0 \in \min(T)$. Suppose that $i < n-1$ and that we already have $t_j$ for each $j \leq i$. Then we choose $t_{i+1}$ as follows. 

If   there is some $t_j \leq t_i$ and there is an immediate successor of $t_j$ which is not in $\{t_k : k \leq i\}$, then let $t_j \leq t_i$ be the $\leq$-largest such $t_j$ and let $t_{i+1}$ be an immediate successor of $t_j$ that is not in $\{t_k : k \leq i\}$. Otherwise, let $t_{i+1} \in \min(T) \backslash \{t_k : k \leq i\}$. 

Let $\tri$ be such that $t_i \tri t_j$ iff $i < j$. One easily verifies that $t_0,t_1, \ldots t_{n-1}$ and $\tri$ are as required. \qed

\bigskip

If $(T,<)$ is a finite forest, then for each $t \in T$, we let $x_t$ be a variable so that $x_s,x_t$ are distinct variables whenever $s,t \in T$ are distinct. 
We say that $\langle \varphi_t(\overline x) : t \in T \rangle$ is a $(T,<)$-{\bf indexed set of formulas} if  the following:
\begin{itemize}

\item $(T,<)$ is a finite forest;
\item each $\varphi_t(\overline x)$ is an ${\mathcal L}(M)$-formula with free variables among \\ $\{x_s : s \in D(t)\}$;

\item if $t \in T$, $|D(t)| = n+1$ and $D(t) = \{t_0,t_1, \ldots, t_{n-1},t\}$, then
$$
\MM \models \forall x_{t_0},x_{t_1}, \ldots, x_{t_{n-1}} \exists x_t \varphi_t(\overline x).
$$
\end{itemize}
If $\langle \varphi_t(\overline x) : t \in T \rangle$ is a $(T,<)$-indexed set of formulas, then  $\langle b_t : t \in T \rangle$ is a {\bf realization of} $\langle \varphi_t(\overline x) : t \in T \rangle$ if 
$$
\MM \models \varphi_t(\overline b)
$$
for each $t \in T$.

\bigskip 

The next lemma generalizes Lemma~3.1. (We easily see that it is indeed  a generalization by considering the empty forest.) 
\bigskip

{\sc Lemma 3.6}: {\em Suppose that  $K',X \subseteq M$ are finite. Let  $\langle \varphi_t(\overline x) : t \in T \rangle$ be a $(T,<)$-indexed set of ${\mathcal L}\big(\acl(K')\big)$-formulas. Then there are a realization $\langle b_t : t \in T \rangle$ and $\NN \preccurlyeq \MM$ such that  whenever $W \subseteq T$ is a twig and $K = K' \cup \{b_t : t \in W\}$, then
\begin{equation} \tag{$***$}
\acl(K  \cup X) \cap N = \acl(K)
\end{equation}
and such that $(\MM,N)$ is arithmetically saturated.}

\bigskip

{\it Proof}. Let $K',X, (T,<)$ and $\langle \varphi_t(\overline x) : t \in T \rangle$ be as given.  By the nullification of $K'$, assume that $K' = \varnothing$. 
Suppose that $|T| = n < \omega$. Then let $\tri$ and $t_0 \tri t_1 \tri \cdots \tri t_{n-1}$ be as in Lemma~3.5. We will obtain three sequences by recursion: $\langle \NN_i : i \leq n \rangle$, 
 $\langle K_i : i \leq n \rangle$ and $\langle d_i : i < n \rangle$.

We start off by applying Lemma~3.1 and letting $\NN_0 \preccurlyeq \MM$ be such that $\acl(X) \cap N_0 = \acl(\varnothing)$ and
$(\MM,N_0)$ is arithmetically saturated. Let $K_0 = \varnothing$. 
Thus,
$$
\acl(K_0 \cup X) \cap N_0 = \acl(K_0).
$$
We continue recursively so that for each $i < n$, we have that:  
\begin{itemize}
\item[(1.i)] $\NN_{i+1} \preccurlyeq \NN_i$ and $(\MM, N_{i+1})$ is arithmetically saturated;
\item[(2.i)] $d_i \in N_i$ and $\MM \models \varphi_{t_i}(\overline b)$, where $b_{t_j} = d_j$ for $t_j \in D(t_i)$;
\item[(3.i)] $K_{i+1} = \{d_i : i \in D(t_i)\}$.
\end{itemize}
Suppose that  $i < n$ and that we already have $\langle \NN_j: j \leq i \rangle$,  
$\langle K_j : j \leq i \rangle$ and $\langle d_j: j < i \rangle$ satisfying (1.j), (2.j) and (3.j) for each  $j<i$. Then choose $d_i \in N_i$, $K_{i+1} = \{d_j : j \leq i\}$ and $\NN_{i+1} \preccurlyeq \NN_i$ such that 
(1.i+1), (2.i) and (3.i+1) are satisfied. 

Thus, we have $\langle \NN_i : i \leq n \rangle$, 
 $\langle K_i : i \leq n \rangle$ and $\langle d_i : i < n \rangle$ such that (1.i), (2.i) and (3.i) are satisfied for all $i < n$. Let $b_{t_i} = d_i$ for all $i < n$, and let $\NN = \NN_n$. 
 The conclusion of Lemma~3.6 easily follows. \qed

\bigskip

Next, we get the corresponding generalization of Corollary~3.2. 

\bigskip

{\sc Corollary 3.7}: {\em Suppose that  $K',X \subseteq M$ are finite. Let  $\langle \varphi_t(\overline x) : t \in T \rangle$ be a $(T,<)$-indexed set of ${\mathcal L}\big(\acl(K')\big)$-formulas. Then there are a realization $\langle b_t : t \in T \rangle$ and $\NN \preccurlyeq \MM$ such that  whenever $W \subseteq T$ is a twig, $K = K' \cup \{b_t : t \in W\}$ and $d \in N$, then
\begin{equation} \tag{$*$}
\acl(K \cup X) \cap \acl(K \cup \{d\}) = \acl(K)
\end{equation}
and such that $(\MM,N)$ is arithmetically saturated.} \qed

\bigskip

Next is  the corresponding generalization of Lemma~2.1. 

\bigskip

{\sc Corollary 3.8}: {\em Suppose that $K',X \subseteq M$ are finite and $a \in M$. Let  $\langle \varphi_t(\overline x) : t \in T \rangle$ be a $(T,<)$-indexed set of ${\mathcal L}\big(\acl(K')\big)$-formulas. Then there are a realization $\langle b_t : t \in T \rangle$ and $d \in M$ such that whenever $W \subseteq T$ is a twig and  $K = K' \cup \{b_t : t \in W\}$, then
\begin{equation} \tag{$*$}
\acl(K \cup  X) \cap \acl(K \cup \{d\}) = \acl( K).
\end{equation}}

\bigskip

{\it Proof}. This follows from Corollary~3.7 in  the same way as Lemma~2.1 follows from Corollary~3.2. \qed
\bigskip

The previous corollary generalized Lemma~2.1. The next lemma is the corresponding generalization of Lemma~2.2. In fact, it is a generalization of an ``over K" version of Lemma~2.2.

\bigskip

{\sc Lemma 3.9}: {\em Suppose that $A',K' \subseteq M$ are finite and $\sigma': A' \into M$ is elementary over $\acl(K')$. Let $\langle \varphi_t(\overline x) : t \in T \rangle$ be a $(T,<)$-indexed set of ${\mathcal L}(K')$-formulas. If $c \in M$, then there are $d \in M$ and a realization 
$\langle b_t : t \in T \rangle$ such that whenever $W \subseteq T$ is a twig, $\sigma \supseteq \sigma'$ and $A = \dom(\sigma)  = A' \cup \{b_t : t \in W\}$, then  $\sigma \cup \{\langle c,d \rangle \}$ is elementary over $K' \cup (\{b_t : t \in W\})$ and  
\begin{equation} \tag{$**$}
\acl(A \cup \{c\}) \cap (\sigma[A] \cup \{d\}) \subseteq \acl(\sigma[A]).
\end{equation}}

\bigskip

{\it Proof}. This follows from Lemma~3.8 in  the same way as Lemma~2.2 follows from Lemma~2.1. \qed

\bigskip

We now come to the crucial definition for proving  the Fix Theorem.

\bigskip

{\sc Definition 3.10}:  We say that $\sigma$ is {\bf superDuby over} $K$ if  every $(T,<)$-indexed set $\langle \varphi_t(\overline x) : t \in T \rangle$ of ${\mathcal L}\big(\acl(K)\big)$-formulas has  a realization  $\langle b_t : t \in T \rangle$ such that whenever $W \subseteq T$ is a twig, then $\sigma$ is Duby over $K \cup \{b_t : t \in W\}$ (as in Definition~2.6). 
If $\sigma$ is superDuby over $\varnothing$,  then we say simply that $\sigma$ is {\bf superDuby}.

\bigskip
In particular, if $\sigma$ is superDuby over $K$, then $\sigma$ is Duby over $K$. 
The following three lemmas, 3.11 -- 3.13,  are used in the proof of the Fix Theorem.

\bigskip

{\sc Lemma 3.11}: {\em The empty function $\varnothing$ is superDuby.}

\bigskip

{\it Proof}. In fact,  for every  finite $K \subseteq M$,
$\varnothing$ is superDuby over $K$. Lemma~3.11 easily follows. \qed

\bigskip

{\sc Lemma 3.12}: {\em Suppose that $\sigma$ is superDuby over $K$ and $a \in M$. Then there is $b \in M$ such that $\tp\big(b / \acl(K)\big) = \tp\big(a / \acl(K)\big)$ and $\sigma$ is superDuby over $K \cup \{b\}$.}

\bigskip

{\it Proof}. Suppose that $\sigma$, $K$ and $a$ are as given. By the nullification of $K$, assume that $K = \varnothing$. Thus, $\sigma$ is superDuby. Consider any nonempty finite forest $(T,<)$, and let $\tri$ be as in Lemma~3.5. Let $t_0 \in T$ be the $\tri$-least element. Let $\varphi(x) \in \tp\big(a / \acl(\varnothing)\big)$. Then let $\langle \varphi_t(\overline x) : t \in T \rangle$ be a $(T,<)$-set 
of ${\mathcal L}\big(\acl(\varnothing)\big)$-formulas such that 
$\varphi_{t_0}(\overline x) = \varphi(x)$. Since $\varphi(x)$ was chosen to be an arbitrary element of $\tp\big(a / \acl(\varnothing)\big)$, by the arithmetic saturation of $\MM$, there is such a realization for which $b_{t_0} \in \tp\big(a / \acl(\varnothing)\big)$. Again, by the arithmetic saturation of $\MM$, there is a single $b_{t_0}$ that works for every such choice of $(T,<)$ and   $\langle \varphi_t(x) : t \in T \rangle$.  Let $b = b_{t_0}$. This $b$ is as in the lemma. \qed

\bigskip

{\sc Lemma 3.13}: {\em  Suppose that  $\sigma$ is superDuby over $K$  and that $c \in M$.}

(a)  {\em There is $d \in M$ such that $\sigma \cup \{\langle c,d \rangle\}$ is superDuby over $K$.} 

(b) {\em There is $d \in M$ such that $\sigma \cup \{\langle d,c \rangle\}$ is superDuby over $K$.}

\bigskip

{\it Proof}. Let $\sigma$, $K$ and $c$ be as given. Just as in Lemma~2.7, it suffices to prove only (a). By the nullification of $K$, assume that $K = \varnothing$. Thus, $\sigma$ is superDuby, and we want $d$ so that $\sigma \cup \{\langle c,d \rangle \}$ is superDuby.

Let   $\langle \varphi_t(\overline x) : t \in T \rangle$ be a $(T,<)$-indexed set of ${\mathcal L}\big(\acl(\varnothing)\big)$-formulas. Then, by Lemma~3.9, there are $d \in M$ and a realization 
$\langle b_t : t \in T\rangle$ as in that lemma. But then, just as in the proof of Lemma~2.4, we get that $d$ is as in (a). Thus, we have that $\sigma \cup \{\langle c,d \rangle\}$ is superDuby, This proves (a) and, therefore, also (b) and Lemma~3.13. \qed

\bigskip

{\it Proof of the Fix Theorem}. (Refer to the  beginning of the Introduction for a statement of the Fix Theorem.) Let $D \subseteq M$ be algebraically closed. If there is a finite $K \subseteq M$ such that $D = \acl(K)$, then the conclusion of the Fix Theorem would follow from Corollary~2.5. Thus, for the remainder of this proof, we  assume that $D$ is infinite. Consequently, $M$ is also infinite. 
Since both $D$ and $M$ are countably infinite, we let $a_0,a_1,a_2, \ldots$   and   $c_0,c_1,c_2, \ldots$  be  nonrepeating enumerations  of $D$ and $M$, respectively.

We will obtain  two infinite sequences: $\langle \sigma_n : n < \omega \rangle$ and $\langle K_n : n < \omega \rangle$. These sequences will be such that whenever $n < \omega$, then:

\begin{itemize}
\item[$(1.n)$] $\sigma_n$ is superDuby over $K_n$;

\item[$(2.n)$] $\sigma_n \subseteq \sigma_{n+1}$ and 
$K_n \subseteq K_{n+1}$.

\end{itemize}

We start by letting $\sigma_0 = K_0 = \varnothing$,  so that by Lemma~3.11, we have that $(1.0)$ is satisfied. We proceed recursively. So, fix some $n< \omega$ and assume that we already have  $\sigma_n$ and $K_n$ satisfying $(1.n)$. Using Lemmas~3.13(a), 3.13(b) and~3.12, we get $d_n,d'_n$ and $b_n \in M$   such that:
\begin{itemize}
\item[$(3.n)$] $\tp\big(b_n / \acl(K_n)\big) = \tp\big(a_n / \acl(K_n)\big)$;
\item[$(4.n)$] $\sigma_n \cup \{\langle c_n,d_n \rangle, \langle d'_n, c_n \rangle\}$ is superDuby over $K_n \cup \{b_n\}$.
\end{itemize}
We then let $\sigma_{n+1} = \sigma_n \cup  \{\langle c_n,d_n \rangle, \langle d'_n, c_n \rangle\}$ and $K_{n+1} = K_n \cup \{b_n\}$.
 It is clear that $(1.n+1)$ and $(2.n)$ are satisfied.

Therefore, by induction, we have these two sequence such that $(1.n)$ and $(2.n)$ are true for each $n < \omega$. 
Let $f = \bigcup_{n<\omega}\sigma_n$. Then, just as in the proof of Duby's Theorem, we see, 
from $(4.n)$ that  $c_n \in \dom(f)$ (by Lemma~3.13(a)) and that $c_n \in \ran(f)$ (by Lemma~3.13(b)), so  that $f \in \aut(\MM)$. Furthermore, each $K_n \subseteq \fix(f)$; in fact,  we have that 
$\fix(f) = \bigcup_{n<\omega}K_n$. Then, $\MM | \fix(f) \cong \MM | D$ as demonstrated by
the function $\{\langle b_n,a_n \rangle : n < \omega\}$. \qed

\bigskip

\bibliographystyle{plain}

\end{document}